\documentclass{article}
\usepackage{amssymb}
\usepackage{amsmath}

\newtheorem{theorem}{Theorem}

\newtheorem{corollary}[theorem]{Corollary}

\newtheorem{definition}[theorem]{Definition}

\newenvironment{proof}[1][Proof]{\noindent\textbf{#1.} }{\ \rule{0.5em}{0.5em}}
\input{tcilatex}
\begin{document}

\title{Some Special Helices in Myller Configuration }
\author{Ak{\i}n Alkan$^1$, Mehmet \"{O}nder$^2$ \\
$^1$Manisa Celal Bayar University, \\
G\"{o}rdes Vocational School, 45750, G\"{o}rdes, Manisa, Turkey.\\
$^2$Delibekirli Village, Tepe Street, 31440, K{\i}r{\i}khan, Hatay, Turkey.\\
E-mails: $^1$akin.alkan@cbu.edu.tr, $^2$mehmetonder197999@gmail.com \\
Orcid Ids: $^1$https://orcid.org/0000-0002-8179-9525, \\
$^2$https://orcid.org/0000-0002-9354-5530}
\maketitle

\begin{abstract}
Some new kinds of special curves called $\overset{\_}{\xi}$-helix, $\overset{%
\_}{\xi}_{1}$-helix, $\overset{\_}{\mu}$-helix, $\overset{\_}{\nu}$-helix
and $W_{k}$-Darboux helices $(k\in \left\{ {n,r,o}\right\} )$ in the Myller
configuration $M(C,\overset{\_}{\xi },\pi )$ are defined and studied. The
necessary and sufficient conditions for these curves are obtained, also the
axes of those helices are given and the relationships between them are
introduced.
\end{abstract}

\textbf{AMS Classsification:} 53A04.

\textbf{Keywords:} Myller configuration; $\overset{\_}{\xi }$-helix; $%
\overset{\_}{\xi }_{1}$ -helix; $W_{k}$-Darboux helices.

\bigskip

\section{Introduction}

\bigskip In the Euclidean 3-space $E^{3}$ the geometric properties of a
curve $C$ are investigated by the aid of othonormal frames defined along the
curve. The well-known one of such frames is the Serret-Frenet frame $\left\{
t,n,b\right\} $, where the unit vector fields $t,n,b$ denote the tangent,
principal normal and binormal of $C$, respectively. Myller considered some
more general frames along a curve $C$. Considering a unit vector $\overset{\_%
}{\xi }$ and a plane $\pi $ along a curve $C$ and calling them as a versor
filed $(C,\overset{\_}{\xi })$ and a plane field $(C,\pi )$ such that~$%
\overset{\_}{\xi }\in \pi ,$ he defined a configuration in$~E^{3}$ called
Myller configuration and denoted by $M(C,\overset{\_}{\xi },\pi )$\cite%
{Miron}. If the planes $\pi ~$are tangent to $C$, then this configuration is
called a tangent Myller configuration and denoted by $M_{t}(C,\overset{\_}{%
\xi },\pi )$. Especially, if the curve $C$ is a surface curve lying on a
surface $S\subset E^{3}$ with arclength parameter $s$, $\overset{\_}{\xi }(s)
$ is the tangent versor field to $S$ along $C$, $\pi (s)$ is the tangent
plane field to $S$ along $C$, then $M_{t}(C,\overset{\_}{\xi },\pi )$ is
called the tangent Myller configuration intrinsic associated to the
geometric objects $S,~C,$ $\overset{\_}{\xi }$. Then, the geometry of the
versor field $(C,\overset{\_}{\xi })$ on a surface $S$ is the geometry of
the associated Myller configurations $M_{t}(C,\overset{\_}{\xi },\pi )$ and
moreover, $M_{t}(C,\overset{\_}{\xi },\pi )$ represents a particular case of 
$M(C,\overset{\_}{\xi },\pi )$. In the special case that tangent Myller
configuration $M_{t}(C,\overset{\_}{\xi },\pi )$ is the associated Myller
configuration to a curve $C$ on a surface $S,$ the classical theory of
surface curves (curves lying on a surface) is obtained.

\bigskip The parallelism of versor field $(C,\overset{\_}{\xi })$ in the
plane field $(C,\pi )$ was studied by Alexandru Myller\cite{Myller}. He
obtained a generalization of parallelism of Levi-Civita on the curved
surfaces. Later, Mayer gave some new invariants for $M(C,\overset{\_}{\xi }%
,\pi )$ and $M_{t}(C,\overset{\_}{\xi },\pi )$\cite{Mayer}. Miron extended
the notion of Myller configuration in Riemannian Geometry\cite{Miron}.
Vaisman considered Myller configuration in the symplectic geometry\cite%
{Vaisman}. Furthermore, a Myller configuration was studied in differenet
spaces \cite{Miron,Constantinescu}. Recently, Macsim, Mihai and Olteanu have
studied rectifying-type curves in a Myller configuration\cite{Macsim}.

In the present paper, we study some special helices in Myller configuration.
We give characterizations for these curves and give the axes. Moreover, we
introduce the relations between these special helices in Myller
configuration $M$.

\section{Preliminaries}

This section gives a brief summary of curves in Myller configuration $M(C,%
\overset{\_}{\xi },\pi ).$ For more detailed information, we refer to \cite%
{Miron}.

\bigskip Let $(C,\overset{\_}{\xi })$ be a versor field and $(C,\pi )$ be a
plane field such that~$\overset{\_}{\xi }\in \pi $ in $E^{3}$. Then, the
pair $((C,\overset{\_}{\xi }),(C,\pi ))$ is called Myller configuration in $%
E^{3}$ and denoted by $M(C,\overset{\_}{\xi },\pi )$ \ or briefly $M$. Let $%
R=(O;\overset{\_}{i}_{1},\overset{\_}{i}_{2},\overset{\_}{i}_{3})$ be an
orthonormal frame. Then, $(C,\overset{\_}{\xi })$ can be given by

\begin{eqnarray*}
\overset{\_}{r} &=&\overset{\_}{r}(s),~\overset{\_}{\xi }=\overset{\_}{\xi }%
(s),~s\in I=(s_{1},s_{2}), \\
\overset{\_}{r}(s) &=&x(s)\overset{\_}{i}_{1}+y(s)\overset{\_}{i}_{2}+z(s)%
\overset{\_}{i}_{3}=\overrightarrow{OP}(s), \\
\overset{\_}{\xi }(s) &=&\overset{\_}{\xi }_{1}(s)\overset{\_}{i}_{1}+%
\overset{\_}{\xi }_{2}(s)\overset{\_}{i}_{2}+\overset{\_}{\xi }_{3}(s)%
\overset{\_}{i}_{3}=\overrightarrow{PQ},
\end{eqnarray*}%
where $s$ is the arclength parameter of the curve $C,$ $\overset{\_}{r}=%
\overset{\_}{r}(s)$ is the position vector of $C$ and $\left\Vert \overset{\_%
}{\xi }(s)\right\Vert ^{2}=\left\langle \overset{\_}{\xi }(s),\overset{\_}{%
\xi }(s)\right\rangle =1.$ Writing $\overset{\_}{\xi }_{1}(s)=\overset{\_}{%
\xi }(s)$ and taking into account $\frac{d\overset{\_}{\xi }_{1}(s)}{ds},$
we define versor field $\overset{\_}{\xi }_{2}(s)$ as follows,%
\begin{equation*}
\frac{d\overset{\_}{\xi }_{1}(s)}{ds}=K_{1}(s)\overset{\_}{\xi }_{2}(s),
\end{equation*}%
where $K_{1}(s)=\left\Vert \frac{d\overset{\_}{\xi }}{ds}\right\Vert $ is
called curvature (or $K_{1}$-curvature) of $(C,\overset{\_}{\xi })$.
Clearly, Since $\overset{\_}{\xi }_{2}(s)$ is orthogonal to $\overset{\_}{%
\xi }_{1}(s)$ and exists when $K_{1}(s)\neq 0$. If we define $\overset{\_}{%
\xi }_{3}(s)=\overset{\_}{\xi }_{1}(s)\times \overset{\_}{\xi }_{2}(s),$ the
frame $R_{F}\left( P(s);\overset{\_}{\xi }_{1}(s),\overset{\_}{\xi }_{2}(s),%
\overset{\_}{\xi }_{3}(s)\right) $ is positively oriented and orthonormal
and called the invariant frame of Frenet-type of the versor field $(C,%
\overset{\_}{\xi })$\cite{Miron}.

The derivative formulas of $R_{F}$ are

\begin{equation*}
\frac{d\overset{\_}{r}(s)}{ds}=\overset{\_}{\alpha }(s)=a_{1}(s)\overset{\_}{%
\xi }_{1}(s)+a_{2}(s)\overset{\_}{\xi }_{2}(s)+a_{3}(s)\overset{\_}{\xi }%
_{3}(s),
\end{equation*}%
with $a_{1}^{2}(s)+a_{2}^{2}(s)+a_{3}^{2}(s)=1$ and

\begin{eqnarray*}
\frac{d\overset{\_}{\xi }_{1}(s)}{ds} &=&K_{1}(s)\overset{\_}{\xi }_{2}(s),
\\
\frac{d\overset{\_}{\xi }_{2}(s)}{ds} &=&-K_{1}(s)\overset{\_}{\xi }%
_{1}(s)+K_{2}(s)\overset{\_}{\xi }_{3}(s), \\
\frac{d\overset{\_}{\xi }_{3}(s)}{ds} &=&-K_{2}(s)\overset{\_}{\xi }_{3}(s),
\end{eqnarray*}%
where $K_{1}(s)>0$ and $K_{2}(s)=\left\langle \frac{d\overset{\_}{\xi }%
_{3}(s)}{ds},\overset{\_}{\xi }_{3}(s)\right\rangle $ is called torsion (or $%
K_{2}$-torsion) of $(C,\overset{\_}{\xi }).$ One can consider the
geometrical interpretation of functions $K_{1}(s)$ and $K_{2}(s)$ as the
same as the curvature and torsion of a curve in $E^{3}$ and the functions $%
a_{1}(s),~a_{2}(s),~a_{3}(s),~K_{1}(s),~K_{2}(s),~(s\in I)~$are invariants
of the versor field $(C,\overset{\_}{\xi }).$ Obviously, if $%
a_{1}(s)=1,~a_{2}(s)=0,~a_{3}(s)=0$, one obtains the Frenet equations of a
curve in $E^{3}$\cite{Miron}.

The following theorem is the fundamental theorem for the versor field $%
(C,\xi )$:

\begin{theorem}
(\cite{Miron})If the functions $K_{1}(s),K_{2}(s),a_{1}(s),a_{2}(s),a_{3}(s),
$ $\left( a_{1}^{2}(s)+a_{2}^{2}(s)+a_{3}^{2}(s)=1\right) $ of class $%
C^{\infty }$ are a priori given, $s\in \lbrack a,b]$, then there exist a
curve $C:[a,b]\rightarrow E^{3}$ with arclength $s$ and a versor field $%
\overset{\_}{\xi }(s)$, $s\in \lbrack a,b]$ such that the functions $a_{i}(s)
$, $(i=1,2,3),$ $K_{1}(s)$ and$~K_{2}(s)$ are the invariants of $(C,\overset{%
\_}{\xi })$. Any two such versor fields $(C,\overset{\_}{\xi })$ differ by a
proper Euclidean motion.
\end{theorem}

We give the definitions of $\overset{\_}{\xi }_{i},(i=1,2,3)$ helices as
follows:\bigskip

\begin{definition}
The curve $C$ with invarinat type Frenet frame $R_{F}\left( P(s);\overset{\_}%
{\xi }_{1}(s),\overset{\_}{\xi }_{2}(s),\overset{\_}{\xi }_{3}(s)\right) $
in $M$ is called $\overset{\_}{\xi }_{i}$-helix if the versor field $\overset%
{\_}{\xi }_{i}$ makes a constant angle with a constant direction, where $%
i=1,2,3$.
\end{definition}

Let $\overset{\_}{v}(s)$ be the unit normal to the oriented plane $\pi .$
Defining $\overset{\_}{\mu }(s)=\overset{\_}{v}(s)\times \overset{\_}{\xi }%
(s),$ positively oriented and orhonormal frame denoted by $R_{D}\left( P(s)%
\text{;}\overset{\_}{\xi }(s),\overset{\_}{\mu }(s),\overset{\_}{v}%
(s)\right) $ and called Darboux frame of the curve $C$ is obtained. This
frame is geometrically associated to Myller configuration $M~$\cite{Miron}.
The following theorem gives the fundamental equations of $M$ .

\begin{theorem}
(\cite{Miron})The moving equations of Darboux frame of $M$ are as follows

\begin{equation*}
\frac{d\overset{\_}{r}}{ds}=\overset{\_}{\alpha }(s)=c_{1}(s)\overset{\_}{%
\xi }(s)+c_{2}(s)\overset{\_}{\mu }(s)+c_{3}(s)\overset{\_}{v}%
(s);~c_{1}^{2}+c_{2}^{2}+c_{3}^{2}=1,
\end{equation*}

and 
\begin{eqnarray*}
\frac{d\overset{\_}{\xi }}{ds} &=&G(s)\overset{\_}{\mu }(s)+K(s)\overset{\_}{%
v}(s), \\
\frac{d\overset{\_}{\mu }}{ds} &=&-G(s)\overset{\_}{\xi }(s)+T(s)\overset{\_}%
{v}(s), \\
\frac{d\overset{\_}{v}}{ds} &=&-K(s)\overset{\_}{\xi }(s)-T(s)\overset{\_}{%
\mu }(s),
\end{eqnarray*}

where $c_{1}(s),~c_{2}(s),~c_{3}(s);~G(s),~K(s)$ and $T(s)$ are invariants
and also uniquely determined. The functions $G(s),~K(s)$ and $T(s)$ are
called the geodesic curvature, the normal curvature and the geodesic torsion
of the versor field $(C,\overset{\_}{\xi )}$, respectively, in $M$.
\end{theorem}

\bigskip

The fundamental theorem for the Darboux frame $R_{D}$ is stated as follows:

\begin{theorem}
(\cite{Miron})Let be a priori given $C^{\infty }$-functions $%
c_{1}(s),~c_{2}(s),~c_{3}(s)$,$~\left[ c_{1}^{2}+c_{2}^{2}+c_{3}^{2}=1\right]
$, $G(s),~K(s),~T(s)$, $s\in \left[ a,b\right] .$ Then, there is a Myller
configuration $M(C,\overset{\_}{\xi },\pi )$ such that the given functions
and parameter $s$ are the invariants and arclength of curve $C,$
respectively. Two such configurations differ by a proper Euclidean motion.
\end{theorem}

\bigskip

The relations between the invarinats of the field $(C,\overset{\_}{\xi })$
and the invarinats of $(C,\overset{\_}{\xi })$ in $M$ are given as follows%
\begin{eqnarray*}
\frac{d\overset{\_}{r}}{ds} &=&\overset{\_}{\alpha }(s)=a_{1}\overset{\_}{%
\xi }_{1}+a_{2}\overset{\_}{\xi }_{2}+a_{3}\overset{\_}{\xi }_{3}=c_{1}%
\overset{\_}{\xi }+c_{2}\overset{\_}{\mu }+c_{3}\overset{\_}{v}, \\
\overset{\_}{\xi }(s) &=&\overset{\_}{\xi }_{1}(s),~\overset{\_}{\xi }%
_{2}(s)=(\sin \psi )\overset{\_}{\mu }(s)+(\cos \psi )\overset{\_}{v}(s),~%
\overset{\_}{\xi }_{3}(s)=-(\cos \psi )\overset{\_}{\mu }(s)+(\sin \psi )%
\overset{\_}{v}(s), \\
c_{1}(s) &=&a_{1}(s),~c_{2}(s)=(\sin \psi )a_{1}(s)-(\cos \psi
)a_{3}(s),~c_{3}(s)=(\cos \psi )a_{2}(s)+(\sin \psi )a_{3}(s),
\end{eqnarray*}

and,%
\begin{equation}
G(s)=(\sin \psi )K_{1}(s),~K(s)=(\cos \psi )K_{1}(s),T(s)=K_{2}(s)+\frac{%
d\psi }{ds},  \label{eqn1}
\end{equation}

where $\psi =\measuredangle (\overset{\_}{\xi }_{2},\overset{\_}{v})$\cite%
{Miron}.

\section{$\protect\overset{\_}{\protect\xi }$-helices in Myller
Configuration $M(C,\protect\overset{\_}{\protect\xi },\protect\pi )$}

In the Euclidean 3-space $E^{3},$ special curves for which their orthonormal
frame vectors make constant angle with some constant directions have an
important role. The well-knowns of such curves studied by Frenet frame are
helices, slant helices and Darboux helices\cite{Struik, Izumiya, Ziplar}.
Moreover, when the curve is a surface curve, then one can also consider
Darboux frame of the curve to study the same special conditions that the
vector fields of the Darboux frame make constant angle with some constant
directions. In this case, surface helices, relatively normal slant helices
and isophote curves are the well-konwn examples of such curves\cite%
{Puig,Dogan,Macit}. In this section, we consider the same conditions in
Myller configuration $M$ and introduced $\overset{\_}{\xi }$-helices in $M$.

\begin{definition}
Let $C$ be a unit speed curve with Darboux frame $R_{D}\left( P(s)\text{;}%
\overset{\_}{\xi }(s),\overset{\_}{\mu }(s),\overset{\_}{v}(s)\right) $ in
Myller configuration $M$. The curve $C$ is called a $\overset{\_}{\xi }$%
-helix in $M$ if the versor field $\overset{\_}{\xi }$ makes a constant
angle with a fixed unit direction $\overset{\_}{d}_{\xi },$ i.e., there
exists a constant angle $\theta $ such that $\left\langle \overset{\_}{\xi },%
\overset{\_}{d}_{\xi }\right\rangle =\cos \theta .$
\end{definition}

\begin{theorem}
\bigskip The curve $C$ with Darboux frame $R_{D}$ and $\left( G,K\right)
\neq (0,0)$ in $M$ is a $\overset{\_}{\xi }$-helix iff the following
function is constant 
\begin{equation}
\sigma _{\xi }=\cot \theta =\mp \frac{K^{2}\left( \frac{G}{K}\right)
^{\prime }-\left( G^{2}+K^{2}\right) T}{\left( G^{2}+K^{2}\right) ^{\frac{3}{%
2}}}.  \label{eqn3}
\end{equation}
\end{theorem}

\begin{proof}
From Definition 5, there exists a unit constant vector $\overset{\_}{d}_{\xi
}$ and a constant function $\theta $ such that $\left\langle \overset{\_}{%
\xi },\overset{\_}{d}_{\xi }\right\rangle =\cos \theta .$ The unit vector $%
\overset{\_}{d}_{\xi }$ can be given in the form 
\begin{equation}
\overset{\_}{d}_{\xi }=(\cos \theta )\overset{\_}{\xi }+x_{2}\overset{\_}{%
\mu }+x_{3}\overset{\_}{v},  \label{3.1}
\end{equation}%
where $x_{2}=x_{2}(s),~x_{3}=x_{3}(s)$ are smooth functions of $s$.
Differentiating (\ref{3.1}) with respect to $s$ gives 
\begin{equation*}
\overset{\_}{d}_{\xi }^{\prime }=\left( -x_{2}G-x_{3}K\right) \overset{\_}{%
\xi }+\left( G\cos \theta +x_{2}^{\prime }-x_{3}T\right) \overset{\_}{\mu }%
+\left( K\cos \theta +x_{2}T+x_{3}^{\prime }\right) \overset{\_}{v},
\end{equation*}%
Since $\overset{\_}{d}_{\xi }$ is constant, we have the system 
\begin{equation}
\left\{ 
\begin{array}{c}
-x_{2}G-x_{3}K=0 \\ 
G\cos \theta +x_{2}^{\prime }-x_{3}T=0 \\ 
K\cos \theta +x_{2}T+x_{3}^{\prime }=0%
\end{array}%
\right.   \label{3.2}
\end{equation}%
From the first equation in (\ref{3.2}), it follows $x_{3}=-x_{2}\frac{G}{K}$%
. Writing that in the second and third equations in (\ref{3.2}) gives
differential equation 
\begin{equation}
x_{2}^{\prime }\left( 1+\left( \frac{G}{K}\right) ^{2}\right) +x_{2}\frac{G}{%
K}\left( \frac{G}{K}\right) ^{\prime }=0.  \label{3.3}
\end{equation}%
The solution of (\ref{3.3}) is $x_{2}=\lambda \frac{K}{\sqrt{G^{2}+K^{2}}}$,
where $\lambda $ is integration constant. Hence, $x_{3}=-\lambda \frac{G}{%
\sqrt{G^{2}+K^{2}}}$. Since $\left\Vert \overset{\_}{d}_{\xi }\right\Vert =1$%
, from (\ref{3.1}) we have $\lambda =\mp \sin \theta .$ Then, (\ref{3.1})
becomes%
\begin{equation}
\overset{\_}{d}_{\xi }=(\cos \theta )\overset{\_}{\xi }\mp \sin \theta \frac{%
K}{\sqrt{G^{2}+K^{2}}}\overset{\_}{\mu }\pm \sin \theta \frac{G}{\sqrt{%
G^{2}+K^{2}}}\overset{\_}{v}.  \label{3.4}
\end{equation}%
By differentiating $\left\langle \overset{\_}{\xi },\overset{\_}{d}_{\xi
}\right\rangle =\cos \theta $ two times, it follows 
\begin{equation*}
-\cos \theta \left( G^{2}+K^{2}\right) \mp \left[ \frac{\left( G^{\prime
}K-GK^{\prime }\right) -\left( G^{2}+K^{2}\right) T}{\sqrt{G^{2}+K^{2}}}%
\right] \sin \theta =0,
\end{equation*}%
which gives 
\begin{equation}
\sigma _{\xi }=\cot \theta =\mp \frac{K^{2}\left( \frac{G}{K}\right)
^{\prime }-\left( G^{2}+K^{2}\right) T}{\left( G^{2}+K^{2}\right) ^{\frac{3}{%
2}}},  \label{3.5}
\end{equation}%
is constant.

Conversely, let the function $\sigma _{\xi }$ given in (\ref{eqn3}) be
constant and $\overset{\_}{d}_{\xi }$ be a unit vector defined by $\overset{%
\_}{d}_{\xi }=\cos \theta \overset{\_}{\xi }\mp \sin \theta \frac{K}{\sqrt{%
G^{2}+K^{2}}}\overset{\_}{\mu }\pm \sin \theta \frac{G}{\sqrt{G^{2}+K^{2}}}%
\overset{\_}{v},$ where $\theta $ is constant. Differentiating last equalty
gives 
\begin{eqnarray*}
\overset{\_}{d}_{\xi }^{\prime } &=&\left( \mp G\sin \theta \left[ \frac{%
-\left( G^{\prime }K-GK^{\prime }\right) +\left( G^{2}+K^{2}\right) T}{%
\left( G^{2}+K^{2}\right) ^{\frac{3}{2}}}\right] +G\cos \theta \right) 
\overset{\_}{\mu } \\
&&+\left( \mp K\sin \theta \left[ \frac{\left( G^{2}+K^{2}\right) T-\left(
G^{\prime }K-GK^{\prime }\right) }{\left( G^{2}+K^{2}\right) ^{\frac{3}{2}}}%
\right] +K\cos \theta \right) \overset{\_}{v}
\end{eqnarray*}%
Now, writing (\ref{eqn3}) in this result, we have $\overset{\_}{d}_{\xi
}^{\prime }=0$, i.e., $\overset{\_}{d}_{\xi }~$ is a constant vector field
and since $\left\langle \overset{\_}{\xi },\overset{\_}{d}_{\xi
}\right\rangle $ is constant, we obtain that $C$ is a $\overset{\_}{\xi }$%
-helix in Myller configuration $M.$
\end{proof}

From Theorem 6, the following corollaries are obtained:

\begin{corollary}
\bigskip The axis of $\overset{\_}{\xi }$-helix $C$ in Myller configuration $%
M$ is given by 
\begin{equation*}
\overset{\_}{d}_{\xi }=(\cos \theta )\overset{\_}{\xi }\mp \sin \theta \frac{%
K}{\sqrt{G^{2}+K^{2}}}\overset{\_}{\mu }\pm \sin \theta \frac{G}{\sqrt{%
G^{2}+K^{2}}}\overset{\_}{v},
\end{equation*}

where $\theta $ is the constant angle defined by $\left\langle \overset{\_}{%
\xi },\overset{\_}{d}_{\xi }\right\rangle =\cos \theta .$
\end{corollary}

\begin{corollary}
i) The curve $C$ with $K=0$ in Myller configuration $M$ is $\overset{\_}{\xi 
}$-helix iff $\sigma _{\xi }=\pm \frac{T}{G}$ is constant.

ii) The curve $C$ with $G=0$ in Myller configuration $M$\ is $\overset{\_}{%
\xi }$-helix iff $\sigma _{\xi }=\pm \frac{T}{K}$ is constant.

iii) The curve $C$ with $T=0$ in Myller configuration $M$ is $\overset{\_}{%
\xi }$-helix iff $\sigma _{\xi }=\mp \frac{G^{\prime }K-GK^{\prime }}{\left(
G^{2}+K^{2}\right) ^{\frac{3}{2}}}$is constant.
\end{corollary}

\begin{theorem}
The curve $C$ in $M$ is a $\overset{\_}{\xi }_{1}$-helix according to the
Frenet-type frame $R_{F}$ iff $\frac{K_{2}}{K_{1}}$ is constant.

\begin{proof}
By taking into account the relations between the invarinats of the field $(C,%
\overset{\_}{\xi })$ and the invarinats of $(C,\overset{\_}{\xi })$ in $M$,
we have $\overset{\_}{\xi }=\overset{\_}{\xi }_{1}.$ Then, it is clear that $%
C$ is a $\overset{\_}{\xi }$-helix if and only if it is a $\overset{\_}{\xi }%
_{1}$-helix. Now, writing (\ref{eqn1}) in (\ref{eqn3}), it follows $\sigma
_{\xi }=\pm \frac{K_{2}}{K_{1}}$, which finishes the proof.
\end{proof}
\end{theorem}

\bigskip

Moreover, from Corollary 7 and Theorem 9, we obtain the following corollary.

\begin{corollary}
\bigskip The axis of a $\overset{\_}{\xi }_{1}$-helix $C$ according to
Frenet-type frame $R_{F}$ is given by $\overset{\_}{d}_{\xi _{1}}=(\cos
\theta )\overset{\_}{\xi }_{1}\pm (\sin \theta )\overset{\_}{\xi }_{3}$,
where $\theta $ is constant.
\end{corollary}

\section{$W_{n}$-Darboux helices in Myller configuration $M(C,%
\protect\overset{\_}{\protect\xi },\protect\pi )$}

Darboux vector which is the angular velocity vector of the Frenet-Serret
frame of a space curve is an important tool to study the differential
geometry of space curves. This vector is directly proportional to angular
momentum that is the reason why it is also called angular momentum vector%
\cite{Stoker}. If we consider another frame along a curve different from the
Frenet frame, the new forms of Darboux vector can be defined. In this
section, we define the Darboux vector of the Darboux frame $R_{D}$ of a
curve $C$ in a Myller configuration $M$. Furthermore, we define some new
forms of Darboux vector and some special Darboux helices in Myller
configuration $M$.

\begin{definition}
Let $C$ be a curve with Darboux frame $R_{D}\left( P(s)\text{;}\overset{\_}{%
\xi }(s),\overset{\_}{\mu }(s),\overset{\_}{v}(s)\right) $ in Myller
configuration $M$. The vector $W=T\overset{\_}{\xi }-K\overset{\_}{\mu }+G%
\overset{\_}{v}$ is called Darboux vector field of the Darboux frame $R_{D}$%
. The vectors%
\begin{equation*}
W_{n}=-K\overset{\_}{\mu }+G\overset{\_}{v},~W_{r}=T\overset{\_}{\xi }+G%
\overset{\_}{v},~W_{o}=T\overset{\_}{\xi }-K\overset{\_}{\mu },
\end{equation*}%
are called the normal-type Darboux vector (or ND-vector), the
rectifying-type Darboux vector (or RD-vector) and the osculating-type
Darboux vector (or OD-vector) of $R_{D}$, respectively.
\end{definition}

\bigskip

Considering versor field $(C,\overset{\_}{W}_{n})$ with $\overset{\_}{W}%
_{n}= $ $\frac{W_{n}}{\left\Vert W_{n}\right\Vert }$, we can give the
followings:

\begin{definition}
Let $C$ be a curve with unit ND-vector field $\overset{\_}{W}_{n}$ in Myller
configuration $M.$The curve $C$ is called $W_{n}$-helix in $M$ if $\overset{%
\_}{W}_{n}$ makes a constant angle with a fixed unit direction $\overset{\_}{%
l}_{n}.$

\begin{theorem}
The curve $C$ with unit ND-vector field $\overset{\_}{W}_{n}$ and $\left(
G,K\right) \neq (0,0)$ in Myller configuration $M$ is a $W_{n}$-helix iff $C$
is a $\overset{\_}{\xi }$-helix in $M.$
\end{theorem}
\end{definition}

\begin{proof}
From Definition 12, we have $\left\langle \overset{\_}{W}_{n},\overset{\_}{l}%
_{n}\right\rangle =\cos \varphi $, where $\varphi $ is the constant angle
between unit versor feleds $\overset{\_}{W}_{n},\overset{\_}{l}_{n}.$
Differentiating that gives%
\begin{equation*}
\frac{K^{2}\left( \frac{G}{K}\right) ^{\prime }-\left( G^{2}+K^{2}\right) T}{%
\left( G^{2}+K^{2}\right) ^{\frac{3}{2}}}\left\langle G\overset{\_}{\mu }+K%
\overset{\_}{v},\overset{\_}{l}_{n}\right\rangle =0.
\end{equation*}%
Since $\overset{\_}{\xi }^{\prime }=G\overset{\_}{\mu }+K\overset{\_}{v},$
last equality becomes $\sigma _{\xi }\left\langle \overset{\_}{\xi }^{\prime
},\overset{\_}{l}_{n}\right\rangle =0.$ If we assume that $\sigma _{\xi }=0$%
, then $\sigma _{\xi }$ is constant, i.e., $C$ is a $\overset{\_}{\xi }$%
-helix in $M.$ Let investigage the case $\sigma _{\xi }\neq 0.$ In this
case, $\left\langle \overset{\_}{\xi }^{\prime },\overset{\_}{l}%
_{n}\right\rangle =0$. On the other hand, by taking into account $%
\left\langle \overset{\_}{W}_{n},\overset{\_}{\xi }^{\prime }\right\rangle
=0 $, we have $\overset{\_}{\xi }^{\prime }=\overset{\_}{W}_{n}\times 
\overset{\_}{l}_{n}.$ Then, it follows that the vectors $\overset{\_}{l}%
_{n},~\overset{\_}{W}_{n}$ and $\overset{\_}{\xi }$ lie on the same plane,
i.e., $\overset{\_}{l}_{n}\in sp\left\{ \overset{\_}{W}_{n},\overset{\_}{\xi 
}\right\} .$ So, we can write $\overset{\_}{l}_{n}=\mp (\sin \varphi )%
\overset{\_}{\xi }+(\cos \varphi )\overset{\_}{W}_{n}$ or, equivalently, 
\begin{equation*}
\overset{\_}{l}_{n}=\mp (\sin \varphi )\overset{\_}{\xi }+\cos \varphi
\left( \frac{-K}{\sqrt{G^{2}+K^{2}}}\mu +\frac{G}{\sqrt{G^{2}+K^{2}}}\overset%
{\_}{v}\right) .
\end{equation*}%
Differentiating $\left\langle \overset{\_}{\xi }^{\prime },\overset{\_}{l}%
_{n}\right\rangle =0,$ gives 
\begin{equation*}
\pm \sin \varphi \left( G^{2}+K^{2}\right) +\cos \varphi \left( \frac{%
-\left( G^{\prime }K-GK^{\prime }\right) +\left( G^{2}+K^{2}\right) T}{%
\left( G^{2}+K^{2}\right) ^{\frac{3}{2}}}\right) =0.
\end{equation*}%
Hence, it follows that 
\begin{equation*}
\cot \varphi =\mp \frac{1}{\frac{K^{2}\left( \frac{G}{K}\right) ^{\prime
}-\left( G^{2}+K^{2}\right) T}{\left( G^{2}+K^{2}\right) ^{\frac{3}{2}}}},
\end{equation*}%
or, equivalently, $\cot \varphi =\mp \frac{1}{\sigma _{\xi }},$which
finishes the proof.
\end{proof}

\bigskip

From Theorem 13, the following corollaries are given.\bigskip

\begin{corollary}
\bigskip The curve $C$ in $M$ is a $\overset{\_}{\xi }_{1}$-helix according
to Frenet-type frame $R_{F}$ iff it is a $W_{n}$-helix.
\end{corollary}

\begin{corollary}
The axis $\overset{\_}{l}_{n}$ of $W_{n}$-helix in $M$ is defined by 
\begin{equation*}
\overset{\_}{l}_{n}=\mp (\sin \varphi )\overset{\_}{\xi }+\cos \varphi
\left( \frac{-K}{\sqrt{G^{2}+K^{2}}}\overset{\_}{\mu }+\frac{G}{\sqrt{%
G^{2}+K^{2}}}\overset{\_}{v}\right) ,
\end{equation*}

where $\varphi $ is the constant angle given by $\left\langle \overset{\_}{W}%
_{n},\overset{\_}{l}_{n}\right\rangle =\cos \varphi .$
\end{corollary}

\section{$\protect\overset{\_}{\protect\mu }$-helices and $\protect\overset{%
\_}{v}$-helices in Myller Configuration $M(C,\protect\overset{\_}{\protect%
\xi },\protect\pi )$}

In this section, we introduce $\overset{\_}{\mu }$-helices, $\overset{\_}{v}$%
-helices, $W_{r}$-Darboux helices and $W_{o}$-Darboux helices in Myller
configuration $M$. The proofs of theorems given below can be made similar to
ones given in previous sections.

\begin{definition}
Let $C$ be a unit speed curve with Darboux frame $R_{D}\left( P(s)\text{;}%
\overset{\_}{\xi }(s),\overset{\_}{\mu }(s),\overset{\_}{v}(s)\right) $ in
Myller configuration $M$. The curve $C$ is called a $\overset{\_}{\mu }$%
-helix in $M$ if the versor field $\overset{\_}{\mu }$ makes a constant
angle with a fixed unit direction $\overset{\_}{d}_{\mu },$ i.e., there
exists a constant angle $\eta $ such that $\left\langle \overset{\_}{\mu },%
\overset{\_}{d}_{\mu }\right\rangle =\cos \eta $.
\end{definition}

\begin{theorem}
The curve $C$ with Darboux frame $R_{D}$ and $\left( G,T\right) \neq (0,0)$
in $M$ is a $\overset{\_}{\mu }$-helix iff the following function is constant%
\begin{equation*}
\sigma _{\mu }=\cot \eta =\mp \frac{G^{2}\left( \frac{T}{G}\right) ^{\prime
}-\left( G^{2}+T^{2}\right) K}{\left( G^{2}+T^{2}\right) ^{\frac{3}{2}}}.
\end{equation*}
\end{theorem}

\begin{corollary}
\bigskip The axis of a $\overset{\_}{\mu }$-helix $C$ in $M$ is given by%
\begin{equation*}
\overset{\_}{d}_{\mu }=\mp \sin \eta \frac{T}{\sqrt{G^{2}+T^{2}}}\overset{\_}%
{\xi +}(\cos \eta )\overset{\_}{\mu }\mp \sin \eta \frac{G}{\sqrt{G^{2}+T^{2}%
}}\overset{\_}{v},
\end{equation*}%
where $\eta $ is the constant angle defined by $\left\langle \overset{\_}{%
\mu },\overset{\_}{d}_{\mu }\right\rangle =\cos \eta .$
\end{corollary}

\begin{corollary}
i) The curve $C$ with $K=0$ in $M$ is $\overset{\_}{\mu }$-helix iff $\sigma
_{\mu }=\mp \frac{T^{\prime }G-TG^{\prime }}{\left( G^{2}+T^{2}\right) ^{%
\frac{3}{2}}}$ is constant.

ii) The curve $C$ with $G=0$ in $M$\ is $\overset{\_}{\mu }$-helix iff $%
\sigma _{\mu }=\pm \frac{K}{T}$ is constant.

iii) The curve $C$ with $T=0$ in $M$ is $\overset{\_}{\mu }$-helix iff $%
\sigma _{\mu }=\pm \frac{K}{G}$ is constant.
\end{corollary}

Considering versor field $(C,\overset{\_}{W}_{r})$ with $\overset{\_}{W}%
_{r}= $ $\frac{W_{r}}{\left\Vert W_{r}\right\Vert }$, we can give the
followings:

\begin{definition}
Let $C$ be a curve with unit RD-vector field $\overset{\_}{W}_{r}$ in Myller
configuration $M.$ The curve $C$ is called a $W_{r}$-helix in $M$\ if $%
\overset{\_}{W}_{r}$ makes a constant angle with a fixed unit direction $%
\overset{\_}{l}_{r}.$
\end{definition}

\begin{theorem}
The curve $C$ with unit RD-vector field $\overset{\_}{W}_{r}$ and $\left(
G,T\right) \neq (0,0)$ in $M$ is a $W_{r}$-helix iff $C$ is a $\overset{\_}{%
\mu }$-helix in $M.$
\end{theorem}

\begin{corollary}
The axis $\overset{\_}{l}_{r}$ of $W_{r}$-helix in $M$ is defined by%
\begin{equation*}
\overset{\_}{l}_{r}=\cos \vartheta \frac{T}{\sqrt{G^{2}+T^{2}}}\overset{\_}{%
\xi }\mp (\sin \vartheta )\overset{\_}{\mu }+\cos \vartheta \frac{G}{\sqrt{%
G^{2}+T^{2}}}\overset{\_}{v},
\end{equation*}%
where $\vartheta $ is the constant angle defined by $\left\langle \overset{\_%
}{W}_{r},\overset{\_}{l}_{r}\right\rangle =\cos \vartheta .$
\end{corollary}

\begin{definition}
Let $C$ be a unit speed curve with Darboux frame $R_{D}$ in Myller
configuration $M$. The curve $C$ is called a $\overset{\_}{v}$-helix in $M$
if the versor field $\overset{\_}{v}$ makes a constant angle with a fixed
unit direction $\overset{\_}{d}_{v},$ i.e., there exists a constant angle $%
\omega $ such that $\left\langle \overset{\_}{v},\overset{\_}{d}%
_{v}\right\rangle =\cos \omega .$
\end{definition}

\begin{theorem}
The curve $C$ with Darboux frame $R_{D}$ and $\left( T,K\right) \neq (0,0)$
in $M$ is a $\overset{\_}{v}$-helix iff the following function is constant%
\begin{equation*}
\sigma _{v}=\cot \omega =\mp \frac{K^{2}\left( \frac{T}{K}\right) ^{\prime
}+\left( T^{2}+K^{2}\right) G}{\left( T^{2}+K^{2}\right) ^{\frac{3}{2}}}.
\end{equation*}
\end{theorem}

\begin{corollary}
The axis of a $\overset{\_}{v}$-helix $C$ in $M$ is given by%
\begin{equation*}
\overset{\_}{d}_{v}=\mp \sin \omega \frac{T}{\sqrt{T^{2}+K^{2}}}\overset{\_}{%
\xi }\pm \sin \omega \frac{K}{\sqrt{T^{2}+K^{2}}}\overset{\_}{\mu }+(\cos
\omega )\overset{\_}{v},
\end{equation*}%
where $\omega $ is the constant angle defined by $\left\langle \overset{\_}{v%
},\overset{\_}{d}_{v}\right\rangle =\cos \omega .$
\end{corollary}

\begin{corollary}
i) The curve $C$ with $K=0$ in $M$ is $\overset{\_}{v}$-helix iff $\sigma
_{v}=\mp \frac{G}{T}$ is constant.

ii) The curve $C$ with $G=0$ in $M$ is $\overset{\_}{v}$-helix iff $\sigma
_{v}=\mp \frac{T^{\prime }K-TK^{\prime }}{\left( T^{2}+K^{2}\right) ^{\frac{3%
}{2}}}$ is constant.

iii) The curve $C$ with $T=0$ in $M$ is $\overset{\_}{v}$-helix iff $\sigma
_{v}=\mp \frac{G}{K}~$is constant.
\end{corollary}

\bigskip

Considering versor field $(C,\overset{\_}{W}_{o})$ with $\overset{\_}{W}%
_{o}= $ $\frac{W_{o}}{\left\Vert W_{o}\right\Vert }$, we can give the
followings:

\begin{definition}
Let $C$ be a curve with unit OD-vector field $\overset{\_}{W}_{o}$ in Myller
configuration $M.$ The curve $C$ is called $W_{o}$-helix in $M$ if $\overset{%
\_}{W}_{o}$ makes a constant angle with a fixed unit direction $\overset{\_}{%
l}_{o}.$
\end{definition}

\begin{theorem}
The curve $C$ with unit OD-vector field $\overset{\_}{W}_{o}$ and $\left(
K,T\right) \neq (0,0)$ in $M$ is an $W_{o}$-helix iff$\ C$ is a $\overset{\_}%
{v}$-helix in $M.$

\begin{corollary}
The axis $\overset{\_}{l}_{o}$ of $W_{o}$-helix in $M$ is defined by%
\begin{equation*}
\overset{\_}{l}_{o}=\cos \varepsilon \frac{T}{\sqrt{T^{2}+K^{2}}}\overset{\_}%
{\xi }-\cos \varepsilon \frac{K}{\sqrt{T^{2}+K^{2}}}\overset{\_}{\mu }\mp
(\sin \varepsilon )\overset{\_}{v},
\end{equation*}%
where $\varepsilon $ is the constant angle defined by $\left\langle \overset{%
\_}{W}_{o},\overset{\_}{l}_{o}\right\rangle =\cos \varepsilon .$
\end{corollary}
\end{theorem}

\section{Conclusions}

Some new types of special curves in Myller configuration $M(C,\overset{\_}{%
\xi },\pi )$ are defined and studied. Definitions and characterizations of $%
\overset{\_}{\xi }$-helices, $\overset{\_}{\mu }$-helices, $\overset{\_}{v}$%
-helices, $W_{k}$-Darboux helices $(k\in \left\{ n,r,o\right\} )$ and $%
\overset{\_}{\xi }_{1}$-helices are introduced. The axes of these helices
are obtained. Also, the relations between these special curves are given.

\textbf{Author Declaration: }There are no known conflicts of interest
associated with this publication and there has been no significant financial
support for this work that could have influenced its outcome.

\textbf{Data Availability Statement:} This manuscript has no associated data.

\end{document}